\documentclass[11pt]{amsart}
\title{Horizontal Linkage of Coherent Functors}
\author{Jeremy Russell}
\address{The College of New Jersey, Ewing, New Jersey}
\email{russelj1@tcnj.edu}
\usepackage{amsmath,amssymb,amsthm, amsfonts, xypic, graphicx}

\usepackage[hidelinks=true]{hyperref}
\usepackage[all]{xy}

\theoremstyle{definition}
\newtheorem{thm}{Theorem}
\newtheorem{prop}[thm]{Proposition}
\newtheorem{lem}[thm]{Lemma}
\newtheorem{cor}[thm]{Corollary}

\newtheorem*{question}{Question}

\newenvironment{pf}{\paragraph{{\sc \dbf{Proof}}}}{\par\vspace{1cm}}
\theoremstyle{definition}
\newtheorem{defn}{\underline{Definition}}

\newtheorem*{thm*}{Theorem}
\newtheorem*{prop*}{Proposition}

\newtheorem*{com}{Comment}
\newtheorem*{lem*}{Lemma}
\newtheorem{ex}{Example}

\renewcommand{\qed}{\blacksquare}

\newcommand{\blank}{\hspace{0.05cm}\underline{\ \ }\hspace{0.1cm} }
\newcommand{\ev}{\textsf{ev}}

\newcommand{\ZZ}{\mathbb{Z}}

\renewcommand{\mod}{\textsf{mod}}
\newcommand{\Mod}{\textsf{Mod}}
\newcommand{\A}{\mathcal{A}}
\newcommand{\B}{\mathcal{B}}
\newcommand{\C}{\mathcal{C}}
\newcommand{\D}{\mathcal{D}}
\newcommand{\E}{\mathcal{E}}

\renewcommand{\P}{\mathcal{P}}

\newcommand{\Y}{\textsf{Y}}

\newcommand{\hr}{\textsf{E}}

\newcommand{\tfp}{\textsf{tfp}}
\newcommand{\fp}{\textsf{fp}}

\newcommand{\dbf}[1]{\textbf{\text{#1}}}

\newcommand{\dia}[1]{\[\xymatrix{#1 }\]}
\usepackage{pgfpages}
\usepackage{tikz}
\usetikzlibrary{matrix,arrows,snakes,backgrounds}

\renewcommand{\:}{\colon}

\newcommand{\To}{\longrightarrow}

\newcommand{\ab}{\textsf{Ab}}

\newcommand{\coker}{\textsf{Coker}}
\renewcommand{\ker}{\textsf{Ker}}

\newcommand{\tor}{\textsf{Tor}}
\newcommand{\ext}{\textsf{Ext}}

\newcommand{\tr}{\textrm{Tr}}

\renewcommand{\u}{\Omega}

\newcommand{\ses}{\textsf{Ses}}

\newcommand{\nat}{\textsf{Nat}}
\renewcommand{\hom}{\textsf{Hom}}

\newcommand{\mcf}{\hom_\Delta}

\begin{document}
\maketitle

\begin{abstract}The satellite endofunctors are used to extend the definition of linkage of ideals to the linkage of totally finitely presented functors.  The new notion for linkage works over a larger class of rings and is consistent with the functorial approach of encoding information about modules into the category of finitely presented functors.  In the process of extending linkage, we recover the Auslander-Gruson-Jensen duality using injective resolutions of finitely presented functors.  Using the satellite endofunctors we give general definitions of derived functors which do not require the existence of projective or injective objects.  A general formula for calculating the defect of a totally finitely presented functor is given.\end{abstract}

\setcounter{tocdepth}{1}
\tableofcontents

\section{Introduction}

The classical notion of horizontal linkage is an ideal theoretic notion and originates in algebraic geometry.  In \cite{linkage}, Martsinkovsky and Strooker provided a new definition of linkage for finitely presented modules over semiperfect rings.  A module $M$ is horizontally linked if $M\cong \u\tr\u\tr M$.  This definition works over a larger class of rings than the classical definition and the original definition can be viewed as a special case of the new definition; however, the operations $\u$ and $\tr$ are not functorial.  

The purpose of this paper is to extend the notion of linkage of modules given by Martsinkovsky and Strooker to a new notion of linkage of finitely presented functors and thereby illustrate how one can recover linkage of modules using functorial methods.  The definition of linkage of functors allows one to predict linkage of modules and works over an even larger class of rings than that considered by Martsinkovsky and Strooker.

The definition is the following.  A functor $F$ is \dbf{totally finitely presented} if there exists finitely presented $R$-modules $X,Y$ and exact sequence $$(Y,\blank)\To (X,\blank)\To F\To 0$$  A totally finitely presented functor $F$ over a coherent ring $R$ is said to be horizontally linked if the counit of adjunction $S^2S_2F\To F$ is an isomorphism.  While this definition is easy to state, the motivation for defining such a functor is not at all straightforward.  In fact certain properties of both the functor category and module category will appear through the process of discovering this definition.  In particular, our search for a more general definition of linkage will bring us to the Auslander-Gruson-Jenson duality.

It turns out that one may arrive at this notion without any motivation from linkage of modules. The satellites form an adjoint pair of endofunctors on the functor category as first shown by Fisher-Palmquist and Newell in \cite{fisheradjoint}.  This naturally leads to the study of those functors $F$ for which the counit of adjunction $S^nS_nF\To F$ is an isomorphism.  Such a question is quite difficult given the size of the functor category and hence a natural approach to this would be to look for a smaller more well behaved functor category.   The category of finitely presented functors would make an excellent starting point; however, even using that category seems too ambitious a start as the left satellite of a finitely presented functor need not be finitely presented even when the ring is coherent.  We are therefore led to the notion of totally finitely presented functors over a coherent ring. 

The paper is organized as follows.  We begin with a brief review of the Yoneda lemma and recall the definition of the category of finitely presented functors.  In that section the definition of finitely presented functors given avoids the larger functor category $(\A,\ab)$ and focuses on building the category $\fp(\A,\ab)$ from the representable functors and the category $\ab$.  This makes everything work in an arbitrary setting.  For a more relaxed definition one may think of a functor $F\:\A\To \ab$ as being finitely presented if there exists an exact sequence of functors $$(Y,\blank)\To (X,\blank)\To F\To 0$$  where exactness is compatible with evaluation at any object in $\A$.  This section appears almost verbatim in \cite{defect} and is included here for the convenience of the reader.

We recall the definitions of the left and right satellites of a functor $F$ in section 3.  These are crucial to our extension of linkage.  Suppose that $\A$ is an abelian category with enough projectives and $\B$ is any abelian category.  The left satellite of any functor $F\:\A\To \B$ is a functor $S_1F\:\A\To \B$ defined as follows.  Given $X\in \A$, take syzygy sequence $$0\To \u X\To P \To X\To 0$$  Then $S_1F(X)$ is completely determined by the exact sequence $$0\To S_1F(X)\To F(\u X)\To F(P)$$  Similarly one can define the dual notion of a right satellite $S^1F$ using cosyzygy sequences.  The basic properties of the satellites needed for this exposition are discussed there.  

The satellite endofunctors form an adjoint pair of endofunctors on the functor category.  Fisher-Palmquist and Newell were the first to show this result and this adjunction is crucial to understanding how to extend the notion of linkage to finitely presented functors.  In addition, the satellite endofunctors produce $\ext$ and $\tor$.  In particular, 
\begin{prop*}For all $n\ge 0$
\begin{enumerate}
\item $S^n(A,\blank)\cong \ext^n(A,\blank)$
\item $S_n\big(\blank\otimes M\big)\cong \tor_n(\blank,M)$
\end{enumerate}
\end{prop*}  

The ideal theoretic notion of linkage of algebraic varieties was extended to a module theoretic notion by Martsinkovsky and Strooker in \cite{linkage}.  In section 4, we recall the definition of horizontal linkage for modules over semiperfect Noetherian rings given there.   Given a finitely presented module $M$, we have a presentation via finitely generated projectives $P_1\To P_0\To M\To 0$.  The transpose of this module is the module $\tr(M)$ appearing in the following exact sequence $$P_0^*\To P_1^*\To \tr(M)\To 0$$ Martsinkovsky and Strooker define a module $M$ over a semiperfect Noetherian ring to be horizontally linked if and only if $M\cong \u\tr\u\tr M$.  Neither the transpose nor the syzygy operation is functorial on $\mod(R)$; however, both are functors on the stable category $\underline{\mod}(R)$. There the objects are the same as $\mod(R)$ and given $X,Y\in \underline{\mod}(R)$, the morphisms are given by $\underline{\hom}(X,Y)=\hom(X,Y)/\P(X,Y)$ where $\P(X,Y)$ is the subgroup of all those morphisms which factor through a projective.   

The objective of the paper is to produce a diagram of functors:   
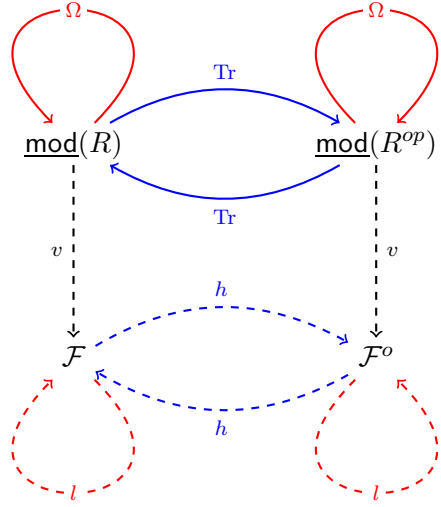
\begin{center}\begin{tikzpicture}[description/.style={fill=white,inner sep=2pt}]  
\matrix (m) [ampersand replacement= \&,matrix of math nodes, row sep=3em, 
column sep=6em,text height=1.5ex,text depth=0.25ex] 
{\underline{\mod}(R)\&\underline{\mod}(R^{op})\\
\&\\
\mathcal{F}\&\mathcal{F}^o\\}; 
\path[->,thick,blue, font=\scriptsize]
(m-1-1) edge[bend left] node[auto]{$\tr$}(m-1-2)
(m-1-2) edge[bend left] node[auto]{$\tr$}(m-1-1);
\path[->,thick, red, font=\scriptsize,  loop/.style={min distance=2cm,in=45, out=135, looseness=13}]
(m-1-2) edge[loop] node[description]{$\Omega$}(m-1-2);
\path[->,thick, red, font=\scriptsize,  loop/.style={min distance=2cm,in=135, out=45, looseness=13}]
(m-1-1) edge[loop] node[description]{$\Omega$}(m-1-1);
\path[->, dashed, thick,blue, font=\scriptsize]
(m-3-1) edge[bend left] node[auto]{$h$}(m-3-2)
(m-3-2) edge[bend left] node[auto]{$h$}(m-3-1);
\path[->, dashed, thick, red, font=\scriptsize,  loop/.style={min distance=2cm,in=-45, out=-135, looseness=13}]
(m-3-2) edge[loop] node[description]{$l$}(m-3-2);
\path[->, dashed,thick, red, font=\scriptsize,  loop/.style={min distance=2cm,in=-135, out=-45, looseness=13}]
(m-3-1) edge[loop] node[description]{$l$}(m-3-1);
\path[->,dashed, thick, font=\scriptsize]
(m-1-1) edge node[left]{$v$} (m-3-1);
\path[->,dashed, thick, font=\scriptsize]
(m-1-2) edge node[auto]{$v$} (m-3-2);
\end{tikzpicture}\end{center} where $\mathcal{F}, \mathcal{F}^o$ are suitable categories of finitely presented functors and $h,l,v$ satisfy certain relations compatible with linkage of modules as defined by Martsinkovsky and Strooker.  We are essentially interested in whether the module theoretic notion of linkage can be recovered as a more general construct in the functor category.

Section 5 is where we get the main results necessary to extend the definition of linkage.   An exact functor $$D_R\:\fp(\Mod(R),\ab)\To\fp(\Mod(R^{op}),\ab)$$ is defined using the fact that $\fp(\Mod(R),\ab)$ has enough injectives and calculating derived functors of the functor $\Y_R$ which is defined by $$\Y_R(F)=\big(F(R),\blank\big).$$   It turns out that all higher derived functors of $\Y_R$ vanish.  We define $$D_R:=L^0\Y_R$$  This functor is shown to be exact; however, it is not a duality.  This leads us to introduce the concept of a totally finitely presented functor.  A functor is \dbf{totally finitely presented}  if there exists $X,Y\in \mod(R)$ and presentation $$(Y,\blank)\To (X,\blank)\To F\To 0$$  Define $\tfp(\Mod(R),\ab)$ to be the full subcategory of $\fp(\Mod(R),\ab)$ consisting of all totally finitely presented subfunctors.  This is an abelian subcategory and the inclusion $$\tfp(\Mod(R),\ab)\To \fp(\Mod(R),\ab)$$ is exact and reflects exact sequences.  

\begin{thm*}The functor  $D_R$ restricts to an exact duality $$D_A\:\tfp(\Mod(R),\ab)\To \tfp(\Mod(R^{op},\ab)$$ on the category of totally finitely presented functors over a coherent ring.  \end{thm*} 

\begin{thm*}The satellite endofunctors $S^n, S_n$ restrict to an adjoint pair of endofunctors on the category of totally finitely presented functors over a coherent ring.  \end{thm*}

The main result of section 6 is that the duality discovered in section 5 is actually the well known Auslander-Gruson-Jensen duality.   As a result, one can recover this duality by exploiting the fact that we have injective resolutions in $\fp(\Mod(R),\ab)$.  

\begin{thm*}The functor $D_A$ is the Auslander-Gruson-Jensen duality.\end{thm*}  While $D_R$ restricts to this duality on the subcategory of totally finitely presented functors, outside that category $D_R$ appears to be different from the functor $D_{aus}$ defined by Auslander: $$D_{aus}F(A):=\nat(F,\blank\otimes A)$$  There is, however, a natural comparison map $D_{aus}\To D_R$.

We also establish the following formula for computing the defect of a totally finitely presented functor over a coherent ring:
\begin{thm*}For any $F\in \tfp(\Mod(R),\ab)$ 
\begin{enumerate}
\item $w\big[D(F)\big]\cong F(R)$
\item $w(F)\cong DF(R)$
\end{enumerate}
\end{thm*}  

The paper concludes by providing a definition of horizontal linkage for totally finitely presented functors over a coherent ring.   A totally finitely presented functor $F$ is horizontally linked if the counit of adjunction $$S^2S_2F\To F$$ is an isomorphism.  The main results of this section are 
\begin{thm*}Let $R$ be a coherent ring.  A module $M$ is horizontally linked in $\underline{\mod}(R)$ if and only if $\ext^1(M,\blank)$ is horizontally linked as a totally finitely presented functor.\end{thm*}

\begin{thm*}All half exact functors of $G$-dimension zero are horizontally linked. \end{thm*}

\begin{thm*}If $F$ is a horizontally linked half exact functor, the $F$ is an extension functor $\ext^1(M,\blank)$.\end{thm*}

The last theorem raises an interesting question.  Are there modules $M$ which are not finitely presented and yet the functor $\ext^1(M,\blank)$ is finitely presented?  In fact, this question is easy to answer because given any finitely presented module $M$, one can take an infinitely generated free module $R^\alpha$ and while $M\oplus R^\alpha$ is not finitely presented, $\ext^1(M\oplus R^\alpha,\blank)\cong \ext^1(M,\blank)$ is finitely presented.  Therefore this question is more accurately a matter of determining whether these are the only examples of such modules.

This paper is one of the main projects appearing in the author's dissertation and many of the results contained herein can be found in \cite{rusthesis}.  The author would like to thank Alex Martsinkovsky for his many helpful suggestions throughout the process of writing it.  In addition, the author would like to thank Ivo Herzog for pointing out that the category of finitely presented functors is a universal construction and informing the author of the deep connections with model theory.  

\section{The Yoneda Lemma and Finitely Presented Functors}  
It is assumed that the reader has some knowledge of abelian categories and functors.  This section is a quick review.  For more details, see \cite{coh} and \cite{mclane}.
Throughout this entire paper:
\begin{enumerate}
\item The word functor will always mean additive functor.    
\item The category of right modules will be denoted by $\Mod(R)$ and the category of finitely presented right modules will be denoted by $\mod(R)$.  
\item Every left module can be viewed as a right module over its opposite ring $R^{op}$.
\item Unless otherwise stated, the category $\A$ will always be assumed to be an abelian category.  This category does not need to be skeletally small.  All statements should work over general abelian categories.
\item The category of abelian groups will be denoted $\ab$.
\end{enumerate}

A functor $F\in(\A, \ab)$ is called \dbf{representable} if it is isomorphic to $\hom_\A(X,\blank)$ for some $X\in \A$.   We will abbreviate the representable functors by $(X,\blank)$.   The most important property of representable functors is the following well known lemma of Yoneda:

\begin{lem}[Yoneda]For any covariant functor $F\:\A\to \ab$ and any $X\in \A$, there is an isomorphism:  $$\nat\big((X,\blank),F\big)\cong F(X)$$  given by $\alpha\mapsto \alpha_X(1_X)$.  The isomorphism is natural in both $F$ and $X$.  \end{lem}

An immediate consequence of the Yoneda lemma is that for any $X,Y\in \A$, $\nat\big((Y,\blank),(X,\blank)\big)\cong (X,Y)$.  Hence all natural transformations between representable functors come from maps between objects in $\A$.   Given a natural transformation between two functors $\alpha:F\To G$, there are functors $\coker(\alpha)$ and $\ker(\alpha)$ which are determined up to isomorphism by their value on any $A\in \A$ by the exact sequence in $\ab$:

$$0\To \ker(\alpha)(A)\To F(A)\overset{\alpha_A}{\To}G(A)\To \coker(\alpha)(A)\To 0$$

A direct consequence of the Yoneda lemma and the fact that $\A$ is abelian is the following:  For any morphism $\alpha\:(Y,\blank)\To (X,\blank)$, the functor $\ker(\alpha)$ is also representable.  Since $\alpha=(f,\blank)$ for some $f\:X\To Y$ there is an exact seqence

Applying $(\blank,A)$ to the exact sequence \begin{center}\begin{tikzpicture}
\matrix(m)[ampersand replacement=\&, matrix of math nodes, row sep=3em, column sep=2.5em, text height=1.5ex, text depth=0.25ex]
{X\&Y\&Z\&0\\}; 
\path[->]
(m-1-1) edge node[auto]{$f$}(m-1-2)
(m-1-2) edge node[auto]{$g$}(m-1-3)
(m-1-3) edge node[auto]{}(m-1-4);
\end{tikzpicture}\end{center} results in the following exact sequence  \begin{center}\begin{tikzpicture}
\matrix(m)[ampersand replacement=\&, matrix of math nodes, row sep=3em, column sep=2.5em, text height=1.5ex, text depth=0.25ex]
{0\&(Z,A)\&(Y,A)\&(X,A)\\}; 
\path[->]
(m-1-1) edge node[auto]{}(m-1-2)
(m-1-2) edge node[auto]{$(g,A)$}(m-1-3)
(m-1-3) edge node[auto]{$(f,A)$}(m-1-4);
\end{tikzpicture}\end{center} Hence, the kernel of any natural transformation between representable functors is itself representable.  Specifically, $$\ker(\alpha)\cong (Z,\blank).$$

\begin{defn}A functor $F\:\A\To \ab$ is \dbf{finitely presented} if there exists $X,Y\in \A$ and $\alpha\:(Y,\blank)\To (X,\blank)$ such that $F\cong \coker(\alpha)$\end{defn}

One easily shows that if $F$ is finitely presented, then the collection of natural transformations $\nat(F,G)$ for any functor $G\:\A\To \ab$ is actually an abelian group.  As such, one may form a category whose objects are the covariant finitely presented functors $F\:\A\To \ab$ and whose morphisms are the natural transformations between two such functors.  This category is denoted by $\fp(\A,\ab)$.  It was first introduced by Auslander and he studied it in great detail in multiple works.  The following is a collection of Auslander's results mainly found in \cite{coh}.

\begin{thm}[Auslander,\cite{coh}]The category $\fp(\A,\ab)$ consisting of all finitely presented functors together with the natural transformations between them satisfies the following properties:
\begin{enumerate}
\item $\fp(\A,\ab)$ is abelian.  The abelian structure is inherited from the category of abelian groups.  More precisely, a sequence of finitely presented functors $$F\To G\To H$$ is exact if and only if for every $A\in \A$ the sequence of abelian groups $$F(A)\To G(A)\To H(A)$$ is exact.
\item $\fp(\A,\ab)$ has enough projectives and they are precisely the representable functors $(X,\blank)$.
\item The functor $\Y\:\A\To \fp(\A,\ab)$ given by $\Y(X):=(X,\blank)$ is a left exact embedding.  This is commonly referred to as the \dbf{Yoneda embedding}.
\item Every finitely presented functor $F\in \fp(\A,\ab)$ has a projective presentation $$0\To (Z,\blank)\To (Y,\blank)\To (X,\blank)\To F\To 0$$
\item The Yoneda embedding $\Y\:\A\To \fp(\A,\ab)$ given by $\Y(X)=(X,\blank)$ is a left exact functor.
\end{enumerate}
\end{thm}

We now give two important examples.

\begin{ex}Suppose that $\A$ has enough projectives.  The functor $\ext^1(A,\blank)$ can be defined without using injective resolutions as follows.  From the short exact syzygy sequence in $\A$: $$0\To \u A\To P\To A\to 0$$ there is a presentation $$0\To(A,\blank)\To (P,\blank)\To (\u A,\blank)\To \ext^1(A,\blank)\To 0$$  As defined $\ext^1(A,\blank)$ is finitely presented.  Moreover, it is easily seen that if $\A$ also has enough injectives, then the abelian group $\ext^1(A,B)$ obtained by using injective resolutions of $B$ is the same as applying the functor $\ext^1(A,\blank)$ as defined above to the object $B$.\end{ex}

In his landmark paper on coherent functors, Auslander showed that given any abelian category $\A$, the category of finitely presented functors $\fp(\A,\ab)$ has some very nice homological properties.  Given a ring $R$, the category $\Mod(R)$ is an abelian category and hence one can look at the category of finitely presented functors $\fp(\Mod(R),\ab)$.   One of the most amazing results appearing in \cite{coh} is a complete description of finitely presented tensor functors $\blank\otimes M$.

\begin{ex}[Auslander, \cite{coh}]Let $R$ denote any ring.  Given a module $M\in \Mod(R^{op})$, the functor $\blank\otimes M\:\Mod(R)\To \ab$ is finitely presented if and only if $M$ is a finitely presented module. \end{ex}

Because this result is sometimes stated in a weaker version, it is worth pointing out the scope of the statement.  If one starts with a finitely presented tensor functor $$(Y,\blank)\To (X,\blank)\To \blank\otimes M\To 0$$ where $X,Y$ are arbitrary not necessarily finitely presented modules, it is possible to show that in fact $M$ is finitely presented.  This means that one can find small modules $A,B\in \mod(R)$ and presentation $$(B,\blank)\To (A,\blank)\To \blank\otimes M\To 0$$

\section{The Satellites}
We begin by recalling the classical definitions of the satellites of any additive functor and listing the main properties. The majority of the proofs can be found in \cite{ce}. Fix abelian category $\C$ with enough projectives and injectives and abelian category $\B$.  Let $F\:\C\to \D$ and $X\in \C$.  Choose syzygy and cosyzygy sequences

\begin{center}\begin{tikzpicture}
\matrix(m)[ampersand replacement=\&, matrix of math nodes, row sep=3em, column sep=2.5em, text height=1.5ex, text depth=0.25ex]
{0\&\Omega X\&P\&X\&0\\
0\&X\&I\&\Sigma X\&0\\};
\path[->]
(m-1-1) edge (m-1-2)
(m-1-2) edge (m-1-3)
(m-1-3) edge (m-1-4)
(m-1-4) edge (m-1-5)
(m-2-1) edge (m-2-2)
(m-2-2) edge (m-2-3)
(m-2-3) edge (m-2-4)
(m-2-4) edge (m-2-5);
\end{tikzpicture}\end{center} The functors $S_1F$ and $S^1F$ are defined on the object $X$ by the following exact sequences:

\begin{center}\begin{tikzpicture}
\matrix(m)[ampersand replacement=\&, matrix of math nodes, row sep=3em, column sep=2.5em, text height=1.5ex, text depth=0.25ex]
{0\&S_1F(X)\&F(\Omega X)\&F(P)\\
F(I)\&F(\Sigma X)\&S^1F(X)\&0\\};
\path[->]
(m-1-1) edge (m-1-2)
(m-1-2) edge (m-1-3)
(m-1-3) edge (m-1-4)
(m-2-1) edge (m-2-2)
(m-2-2) edge (m-2-3)
(m-2-3) edge (m-2-4);
\end{tikzpicture}\end{center} 
These assignments are independent, up to isomorphism, of any choices involved and in fact these sequences completely determine $S_1F$ and $S^1F$ as functors.  Moreover, $S_1$ and $S^1$ are functorial in $F$.  The satellites can be interpreted as homology of a certain complex.  To see this, note that from the above syzygy sequence we have a complex $$0\To \u X\To P\To 0$$ Applying $F$ to this complex gives the complex  $$0\To F(\u X)\To F(P)\To 0$$ at $F(\u X)$.   $S_1F(X)$ is the homology of this complex at $F(\u X)$. Different choices of syzygy sequences will result in homotopic complexes and $F$ will preserve this homotopy. Thus different choices  result in isomorphic homologies. A similar approach can be taken for realizing $S^1F$ as homology.

By iteration, one defines satellites $S^n=S^1S^{n-1}$ and $S_n=S_1S_{n-1}$ for $n>0$.  Notationally it is sometimes more convenient to write $S^{-n}$ in place of $S_n$ for $n>0$.  In this notation the satellites of any functor $F$ form what is known as a connected sequence of functors $(S^nF)^{n\in \ZZ}$ where $S^0F=F$.  This is a sequence of additive covariant functors satisfying the the following property.

\begin{prop}For any morphism of short exact sequences in $\C$:

\begin{center}\begin{tikzpicture}
\matrix(m)[ampersand replacement=\&, matrix of math nodes, row sep=3em, column sep=2.5em, text height=1.5ex, text depth=0.25ex]
{0\&A\&B\&C\&0\\
0\&X\&Y\&Z\&0\\};
\path[->]
(m-1-1) edge (m-1-2)
(m-2-1) edge (m-2-2)
(m-1-2) edge (m-1-3)
edge (m-2-2)
(m-2-2) edge (m-2-3)
(m-1-3) edge (m-1-4)
edge (m-2-3)
(m-2-3) edge (m-2-4)
(m-1-4) edge (m-1-5)
edge (m-2-4)
(m-2-4) edge (m-2-5);
\end{tikzpicture}\end{center} there is a commutative diagram whose rows are complexes:

\begin{center}\begin{tikzpicture}
\matrix(m)[ampersand replacement=\&, matrix of math nodes, row sep=3em, column sep=2.5em, text height=1.5ex, text depth=0.25ex]
{\cdots\&S^{n-1}F(C)\&S^nF(A)\&S^nF(B)\&S^nF(C)\&S^{n+1}F(A)\&\cdots\\
\cdots\&S^{n-1}F(C)\&S^nF(X)\&S^nF(Y)\&S^nF(Z)\&S^{n+1}F(X)\&\cdots\\};
\path[->]
(m-1-1) edge (m-1-2)
(m-2-1) edge (m-2-2)
(m-1-2) edge (m-1-3)
edge (m-2-2)
(m-2-2) edge (m-2-3)
(m-1-3) edge (m-1-4)
edge (m-2-3)
(m-2-3) edge (m-2-4)
(m-1-4) edge (m-1-5)
edge (m-2-4)
(m-2-4) edge (m-2-5)
(m-1-5) edge (m-1-6)
edge (m-2-5)
(m-2-5) edge (m-2-6)
(m-1-6) edge (m-2-6)
edge (m-1-7)
(m-2-6) edge (m-2-7);
\end{tikzpicture}\end{center}  Moreover, the rows are exact if and only if $F$ is half exact.\end{prop}

We will not need the notion of a connected sequence of functors for the purposes of this paper and so we will not go into any more detail concerning these objects; however, it should be pointed out that most of the arguments one uses when studying the satellites involve connected sequences of functors as well as local data.  In fact, anyone that has dealt with $\ext$ and $\tor$ has already seen connected sequences of functors.

\begin{prop}[Cartan-Eilenberg \cite{ce}]Let $\C$ be an abelian category with enough projectives and $A\in \C$. The satellite sequence of $\hom(A,\blank)$ is $\ext(A,\blank)$.  That is, given any $n\in \ZZ$
$$S^n(A,\blank)\cong \ext^n(A,\blank)$$ where $\ext^{n}(A,\blank)=0$ for any $n<0$.
\end{prop}

\begin{prop}[Cartan-Eilenberg \cite{ce}]For any $R$-module, the satellite sequence of the tensor functor $\blank\otimes M$ is $\tor(\blank,M)$.  That is, given any $n\in \ZZ$, $$S_n(\blank\otimes M)\cong \tor_n(\blank,M)$$\end{prop}

In \cite{fisheradjoint}, Fisher-Palmquist and Newell established that the left and right satellites are actually adjoint pairs of endofunctors on the functor category.  In that setting, the source category is assumed only to be skeletally small in order to ensure that no set theoretic issues arise.  The reason for this restriction is to ensure that the collection of natural transformations between two functors is actually an abelian group.   In order to avoid any issues at this point, we state the result of Fisher-Palmquist and Newell under the needed assumptions.  We will later be able to relax the assumptions somewhat.

\begin{thm}[Fisher-Palmquist, Newell, \cite{fisheradjoint}]If $\C$ is a skeletally small abelian category with enough projectives and enough injectives and $\D$ is an abelian category, then for any $n\ge 0$, the satellite functors $(S^n,S_n)$ form an adjoint pair of endofunctors on the functor category $(\C,\D)$.  More precisely, for any functors $F,G\:\C\To \D$ $$\nat(S^nF,G)\cong \nat(F,S_nG)$$\end{thm}  


For any additive category $\C$ with enough projectives, the \dbf{stable category} $\underline{\C}$ is defined as follows.  The objects of $\underline{\C}$ are precisely the objects of $\C$.  Given $X,Y\in\underline{\C}$, define $\hom_{\underline{\C}}(X,Y)=\hom_{\C}(X,Y)/\P(X,Y)$ where $\P(X,Y)$ is the subgroup of $\hom(X,Y)$ consisting of all morphisms that factor through a projective.  One writes $\underline{\hom}(X,Y)$ in place of $\hom_{\underline{\C}}(X,Y)$.

The stable category was studied from the functorial point of view in \cite{stab}.   The methodology employed there is explained as follows. If $\C$ is abelian with enough projectives and injectives, then for any additive functor $F\:\C\To \D$, there are exact sequences of functors  \begin{center}\begin{tikzpicture}
\matrix(m)[ampersand replacement=\&, matrix of math nodes, row sep=3em, column sep=2.5em, text height=1.5ex, text depth=0.25ex]
{L_0F\&F\&\overline{F}\&0\\}; 
\path[->]
(m-1-1) edge node[auto]{}(m-1-2)
(m-1-2) edge node[auto]{}(m-1-3)
(m-1-3) edge node[auto]{}(m-1-4);
\end{tikzpicture}\end{center}   \begin{center}\begin{tikzpicture}
\matrix(m)[ampersand replacement=\&, matrix of math nodes, row sep=3em, column sep=2.5em, text height=1.5ex, text depth=0.25ex]
{0\&\underline{F}\&F\&R^0F\\}; 
\path[->]
(m-1-1) edge node[auto]{}(m-1-2)
(m-1-2) edge node[auto]{}(m-1-3)
(m-1-3) edge node[auto]{}(m-1-4);
\end{tikzpicture}\end{center}   where $\underline{F}$ is called the \dbf{projective stabilization} of $F$ and $\overline{F}$ is called the \dbf{injective stabilization} of $F$.   

The relationship between the stable category of the projective stabilization can be realized by the following.

\begin{ex}Let $X\in \C$.  Set $F=(X,\blank)$ in order to make the notation somewhat easier to understand.  The emphasis here is on ``somewhat'' because this choice might also have the complete opposite effect.  There exists an exact sequence of functors $$L_0(X,\blank)\To (X,\blank)\To \underline{F}\To 0$$  Let $Y\in \C$.  Take any projectives $P,Q$ and any exact sequence $$Q\To P\overset{f}{\To} Y\To 0$$ by definition of $L_0$ and the induced map $L_0F\To F$, we have a commutative diagram with exact rows
\dia{(X,Q)\ar[d]\ar@{=}[r]&(X,Q)\ar[d]\\
(X,P)\ar[d]_{c}\ar@{=}[r]&(X,P)\ar[d]^{(X,f)}\\
L_0F(Y)\ar[d]\ar[r]^{g}&(X,Y)\ar[r]&\underline{F}(Y)\ar[r]&0\\
0} In this diagram the left most column is exact.  Since $\underline{F}$ is $(X,Y)$ modulo the image of $g$ and $gc=(X,f)$, the fact that $c$ is an epimorphism forces the image of $g$ to coincide with the image of $(X,f)$.   Hence the image of $(X,f)$ is the subgroup of all morphisms $h\:X\To Y$ that factor through $P$ via $f$.  It is easily seen that a morphism $h\:X\to Y$ satisfies this condition if and only if it factors through a projective, that is if and only if $h\in \P(X,Y)$.  As a result, $$\underline{F}=\hom_{\underline{\C}}(X,\blank)=\underline{\hom}(X,\blank)$$\end{ex}

One of the major results from \cite{stab} is that the satellites determine the projective and injective stabilization of half exact functors.  

\begin{prop}[Auslander-Bridger, \cite{stab}] For any $X\in \C$ and any half exact functor $F\:\C\to \D$
\begin{enumerate}
\item $\overline{F}\cong S^1S_1(F)$
\item $\underline{F}\cong S_1S^1(F)$ 
\item $S_1S^1\hom(X,\blank)\cong \underline{\hom}(X,\blank)$
\end{enumerate}
\end{prop}

\begin{prop}[Auslander-Bridger, \cite{stab}]Suppose that $\A$ and $\B$ are abelian and that $\A$ has enough projectives and injectives.  For any functor $F\:\A\To \B$.
\begin{enumerate}
\item $S_1\overline{F}\cong S_1 F$
\item $S^1\underline{F}\cong S^1 F$
\item If $F$ vanishes on injectives, then $F=\overline{F}$.
\item $\overline{F}$ vanishes on injectives and is hence its own injective stablilization.
\item If $F$ vanishes on projectives, then $F=\underline{F}$.
\item $\underline{F}$ vanishes on projectives and is hence its own projective 
stabilization.
\end{enumerate}
\end{prop}

Moving from $\C$ to the stable category $\underline{\C}$ is somewhat complicated in the sense that morphisms become more difficult to understand.  This is because we identify any two morphisms whose difference factors through a projective in the source category.  One of the main tools that we have at our disposal is the Hilton-Rees embedding.  Notationally, we have decided to denote this contravariant functor by $\hr$ which does not seem to be standard.  

\begin{thm}[Hilton-Rees, \cite{hr}]The functor $$\hr\:\underline{\C}\to\fp(\C,\ab)$$ given by $\hr(A):=\ext^1(A,\blank)$ is fully faithful and hence is a contravariant embedding. \end{thm}

%

\section{Linkage of Modules}We are finally ready to begin our discussion of horizontal linkage.  Horizontal linkage as defined in algebraic geometry is an ideal theoretic notion. In \cite{linkage}, Martsinkovsky and Strooker showed how this notion can be extended to finitely generated modules over semiperfect Noetherian rings.  This will be the first step towards extending linkage to finitely presented functors.   In order to state the definition we need the well known \dbf{transpose} operation.  For any ring $R$ and any module $M\in \mod(R)$, apply the functor $(\blank)^*=\hom(\blank,R)$ to any presentation of $M$:
\begin{center}\begin{tikzpicture}
\matrix(m)[ampersand replacement=\&, matrix of math nodes, row sep=3em, column sep=2.5em, text height=1.5ex, text depth=0.25ex]
{P_1\&P_0\&X\&0\\};
\path[->]
(m-1-1) edge (m-1-2)
(m-1-2) edge (m-1-3)
(m-1-3) edge (m-1-4);
\end{tikzpicture}\end{center} yielding exact sequence
\begin{center}\begin{tikzpicture}
\matrix(m)[ampersand replacement=\&, matrix of math nodes, row sep=3em, column sep=2.5em, text height=1.5ex, text depth=0.25ex]
{0\&X^*\&(P_0)^*\&(P_1)^*\&\tr(M)\&0\\};
\path[->]
(m-1-1) edge (m-1-2)
(m-1-2) edge (m-1-3)
(m-1-3) edge (m-1-4)
(m-1-4) edge (m-1-5)
(m-1-5) edge (m-1-6);
\end{tikzpicture}\end{center} The module $\tr(M)$ is called the transpose of $M$.  It  depends on the chosen presentation of $M$ in general.  We now recall the definition of linkage given by Martsinkovksy and Strooker.  Let $R$ be a semiperfect Noetherian ring.  A finitely generated module $M$ is \dbf{horizontally linked} if $M\cong \u \tr\u\tr(M)$.

One of the main ideas behind the functorial approach is that one can study the module category by embedding it into the category of finitely presented functors. The definition of linkage given by Martsinkovsky-Strooker is stated in terms of $\u$ and $\tr$ which are module theoretic operations.  The main goal of this paper is to answer the following question:

\begin{question}Can the notion of linkage on $\mod(R)$ be extended to a notion of linkage of finitely presented functors?\end{question}

There are two important issues that must be addressed:
\begin{enumerate}
\item The operations $\u$ and $\tr$ are not functorial constructions.
\item It is not clear what extension of linkage means.
\end{enumerate}

Martsinkovksy made the fundamental observation that the first issue can be solved by passing to the stable module category.  The objects of $\underline{\mod}(R),\ab)$ are finitely presented $R$-modules.  Given $X,Y\in\underline{\mod}(R)$, $$\hom_{\underline{\mod}(R)}(X,Y):=\hom_{\mod(R)}(X,Y)/\P(X,Y)$$ where $\P(X,Y)$ is the subgroup of $\hom(X,Y)$ consisting of those morphisms that factor through projectives..  It is well known that the operations $\u$ and $\tr$ on $\mod(R)$ induce functors $\u$ and $\tr$ on the stable category.  Moreover, $(\tr\u\tr,\u)$ form an adjoint pair of endofunctors.  This means that there is a unit of adjunction \dia{1_{\underline{\mod}(R)}\ar[r]& \u\tr\u\tr}  This allows one to restate the definition of linkage for stable modules:

\begin{defn}Let $R$ be a Noetherian ring.  A module $M\in \underline{\mod}(R)$ is linked if $u_M$ is an isomorphism: \dia{M\ar[r]^{u_M\qquad}_{\cong\qquad}&\u\tr\u\tr(M)}\end{defn}
 With this in mind, we would like to produce 2 suitable categories of finitely presented functors $\mathcal{F},\mathcal{F}^o$ and a diagram of functors:
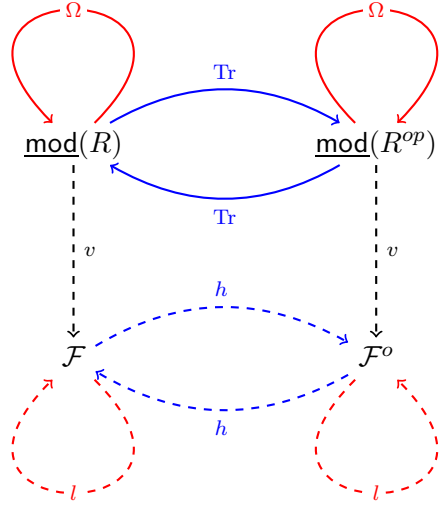
\begin{center}\begin{tikzpicture}[description/.style={fill=white,inner sep=2pt}]  
\matrix (m) [ampersand replacement= \&,matrix of math nodes, row sep=3em, 
column sep=6em,text height=1.5ex,text depth=0.25ex] 
{\underline{\mod}(R)\&\underline{\mod}(R^{op})\\
\&\\
\mathcal{F}\&\mathcal{F}^o\\}; 
\path[->,thick,blue, font=\scriptsize]
(m-1-1) edge[bend left] node[auto]{$\tr$}(m-1-2)
(m-1-2) edge[bend left] node[auto]{$\tr$}(m-1-1);
\path[->,thick, red, font=\scriptsize,  loop/.style={min distance=2cm,in=45, out=135, looseness=13}]
(m-1-2) edge[loop] node[description]{$\Omega$}(m-1-2);
\path[->,thick, red, font=\scriptsize,  loop/.style={min distance=2cm,in=135, out=45, looseness=13}]
(m-1-1) edge[loop] node[description]{$\Omega$}(m-1-1);
\path[->, dashed, thick,blue, font=\scriptsize]
(m-3-1) edge[bend left] node[auto]{$h$}(m-3-2)
(m-3-2) edge[bend left] node[auto]{$h$}(m-3-1);
\path[->, dashed, thick, red, font=\scriptsize,  loop/.style={min distance=2cm,in=-45, out=-135, looseness=13}]
(m-3-2) edge[loop] node[description]{$l$}(m-3-2);
\path[->, dashed,thick, red, font=\scriptsize,  loop/.style={min distance=2cm,in=-135, out=-45, looseness=13}]
(m-3-1) edge[loop] node[description]{$l$}(m-3-1);
\path[->,dashed, thick, font=\scriptsize]
(m-1-1) edge node[auto]{$v$} (m-3-1);
\path[->,dashed, thick, font=\scriptsize]
(m-1-2) edge node[auto]{$v$} (m-3-2);
\end{tikzpicture}\end{center}

such that 
\begin{enumerate}
\item The vertical functors $v$ are contravariant embeddings.
\item The following commutativity relations are satisfied:
\begin{enumerate}
\item $hv=v\tr$
\item $l^khv=v\u^k\tr$
\end{enumerate}
\item There is a counit of adjunction \dia{lhlh\ar[r]^c& 1}
\end{enumerate}

If such a diagram exists, we will define a functor $F$ to be linked if $c_F$ is an isomorphism.  We will then further require that $M$ is linked in $\underline{\mod}(R)$ if and only if $v(M)$ is linked in $\mathcal{F}$.  We will consider such a situation a satisfactory way of extending the definition of linkage.

We now return to the satellites.  The functor $S^1$ is right exact.  Let $F\in \fp(\Mod(R),\ab)$.  Applying the right exact functor $S^1$ to the presentation \dia{(Y,\blank)\ar[r]&(X,\blank)\ar[r]&F\ar[r]&0} produces exact sequence \dia{\ext^1(Y,\blank)\ar[r]&\ext^1(X,\blank)\ar[r]&S^1F\ar[r]&0}  Since $\Mod(R)$ is has enough projectives, the functors $\ext^1(Y,\blank)$ and $\ext^1(X,\blank)$ are finitely presented.  Therefore the functor $S^1F$ is finitely presented.  This means that $S^1\:(\Mod(R),\ab)\to (\Mod(R),\ab)$ restricts to a functor $S^1\:\fp(\Mod(R),\ab)\to \fp(\Mod(R),\ab)$.  Moreover, $$S^1\ext^1(M,\blank)\cong \ext^1(\u M,\blank)$$  Combining this with the Hilton-Rees embedding $\hr$, we start to suspect that we should choose $v=\hr$ and $l=S^1$.  That is, it appears that if we use the Hilton-Rees embedding, we already have an analog for $\u$ on the functor category; however, we are still missing the analog of the transpose.  Moreover, the functor $S^1$ is well behaved on $\fp(\Mod(R),\ab)$ and well behaved on extension functors; however, it is not clear what choice we should make for $\mathcal{F}$ and $\mathcal{F}^o$.  We also cannot ignore the fact that $\tr\:\underline{\mod}(R)\To \underline{\mod}(R^{op})$ is a duality on the stable module category.  Ideally, our new definition of linkage should involve categories $\mathcal{F},\mathcal{F}^o$ that are also dual.

\section{Duality and Large Modules}

%
%

Let $R$ be a ring.  The category $\Mod(R)$ is abelian with enough projectives.  Therefore, as seen in \cite{defect}, the category of finitely presented functors $\fp(\Mod(R),\ab)$ has enough injectives and in fact they are completely classifed as quotients of natural transformations between exact functors.   There is an obvious contravariant functor $$\Y_R\:\fp(\Mod(R),\ab)\to \fp(\Mod(R^{op}),\ab)$$ defined by $$\Y_R(F):=(F(R),\blank)$$  Because $\fp(\Mod(R),\ab)$ has enough injectives, one may calculate the left derived functors $L^n(\Y_R)$.  It is easily seen that for all $n\ge 1$, $L^n(\Y_R)=0$.  The only surviving derived functor is $L^0\Y_R$.  
\begin{defn}The functor $D_R\:\fp(\Mod(R),\ab)\to \fp(\Mod(R^{op}),\ab)$ is defined by $$D_R:=L^0\Y_R$$  \end{defn}

\begin{lem}Suppose that $\C$ is an abelian category with enough injectives and $\D$ is an abelian category.  In addition, assume that every object in $\C$ has injective dimension at most 2.  Then for any left exact contravariant functor $S\:\C\to \D$, the zeroth derived functor $L^0S$ is exact.\end{lem}

\begin{pf}Take any exact sequence $0\to X\to Y\to Z\to 0$ in $\C$.  There is a commutative diagram with exact rows and columns: \dia{&0\ar[d]&0\ar[d]&0\ar[d]&\\
0\ar[r]&X\ar[r]\ar[d]&Y\ar[d]\ar[r]&Z\ar[d]\ar[r]&0\\
0\ar[r]&I^0\ar[d]\ar[r]&J^0\ar[d]\ar[r]&K^0\ar[d]\ar[r]&0\\
0\ar[r]&I^1\ar[d]\ar[r]&J^1\ar[d]\ar[r]&K^1\ar[d]\ar[r]&0\\
0\ar[r]&I^2\ar[d]\ar[r]&J^2\ar[d]\ar[r]&K^2\ar[d]\ar[r]&0\\
&0&0&0&} where $I^p,J^p,K^p$ are injectives.  From this diagram, the fact that $S$ is left exact, and the fact that the rows consisting of injectives split, one easily recovers the following commutative diagram with exact rows and columns:\dia{&0\ar[d]&0\ar[d]&0\ar[d]&\\
0\ar[r]&S(K^2)\ar[r]\ar[d]&S(J^2)\ar[d]\ar[r]&S(I^2)\ar[d]\ar[r]&0\\
0\ar[r]&S(K^1)\ar[d]\ar[r]&S(J^1)\ar[d]\ar[r]&S(I^1)\ar[d]\ar[r]&0\\
0\ar[r]&S(K^0)\ar[d]\ar[r]&S(J^0)\ar[d]\ar[r]&S(I^0)\ar[d]\ar[r]&0\\
&L^0S(Z)\ar[d]\ar[r]&L^0S(Y)\ar[d]\ar[r]&L^0S(X)\ar[d]\ar[r]&0\\
&0&0&0&}  Applying the snake lemma yields exact sequence $$0\to L^0S(Z)\to L^0S(Y)\to L^0S(X)\to 0$$ completing the proof.  $\qed$ 
\end{pf}

\begin{prop}The functor $D_R$ is exact.\end{prop}

\begin{pf}The functor $\Y_R$ is the composition of evaluation at the ring $\ev_R$ followed by the Yoneda embedding $\Y$.  Since evaluation is exact and the Yoneda embedding is left exact, $\Y_R$ is left exact.  Hence $D_R=L^0\Y_R$ is exact by the preceding lemma. $\qed$\end{pf}

A ring is called \dbf{coherent} if the categories $\mod(R)$ and $\mod(R^{op})$ are abelian.  For the remainder of this paper, we will assume that any ring $R$ is coherent.  The versatility of looking at functors $F\in \fp(\Mod(R),\ab)$ instead of restricting to $\fp(\mod(R), \ab)$ is that $\Mod(R)$ is abelain, has enough projectives and injectives and hence satellite functors can be defined using the classical definitions. We will be concerned with finitely presented functors that actually arise from representable functors $(X,\blank)$ for which $X$ itself is a finitely presented module.   We will refer to these functors as totally finitely presented.     

\begin{defn}Let $\tfp(\Mod(R),\ab)$ represent the full subcategory of $\fp(\Mod(R),\ab)$ consisting of all finitely presented functors $F\in \fp(\Mod(R),\ab)$ such that there exists $X,Y\in \mod(R)$ and presentation $$(Y,\blank)\To (X,\blank)\To F\To 0$$ A \dbf{totally finitely presented functor} is any functor $F\in\tfp(\Mod(R),\ab)$.   \end{defn}

\begin{thm}For a coherent ring $R$, the full subcategory of totally finitely presented functors $\tfp(\Mod(R),\ab)$ is abelian and the inclusion $$\tfp(\Mod(R),\ab)\To \fp(\Mod(R),\ab)$$ is exact and reflects exact sequences.  In particular, the notion of exactness in $\tfp(\Mod(R),\ab)$ is compatible with evaluation.  This category is equivalent to $\fp(\mod(R),\ab)$.  \end{thm}   

\begin{pf}Let $\P$ be the full subcategory of $\fp(\Mod(R),\ab)$ consisting of all representable functors $(X,\blank)$ for which $X$ is finitely presented as a module.  This category consists of projectives in $\fp(\Mod(R),\ab)$ and since $R$ is coherent it is closed under kernels and finite sums.  It follows from a result of Auslander appearing in \cite{coh} that the category of all functors with presentations from $\P$ is an abelian subcategory and it is easily seen that the inclusion is exact and reflects exact sequences.   

The equivalence of $\tfp(\Mod(R),\ab)$ and $\fp(\mod(R),\ab)$ is given by the restriction functor $\tfp(\Mod(R),\ab)\to \fp(\mod(R),\ab)$.  The fact that all functors $F\in \tfp(\Mod(R),\ab)$ commute with direct limits and hence are completely determined by  $\mod(R)$ implies that the restriction is an equivalence.  $\qed$\end{pf} 

In \cite{coh}, Auslander establishes that the functor $\blank\otimes M\: \Mod(R)\To\ab$ is finitely presented if and only if $M\in\mod(R^{op})$.  Notice that since $M$ is finitely presented, one can easily verify that $\blank\otimes M$ is an object of $\tfp(\Mod(R),\ab)$.  In other words $\blank\otimes M$ is finitely presented if and only if it is totally finitely presented.  In fact, for a coherent ring we have the following well known result:

\begin{lem}For a coherent ring $R$, if $M$ is a finitely presented module then for all $n\ge 0$, the functor $\tor_n(\blank, M)$ is totally finitely presented.\end{lem}

We now compute $D_R$ for multiple functors.

\begin{lem}Suppose that $X\in \mod(R^{op})$.  Then $D_R(X,\blank)\cong \blank\otimes X$. Here the functor $D_R(X,\blank)$ is being viewed as an object in $\tfp(\Mod(R),\ab)$.\end{lem}

\begin{pf}From the presentation of $X$: 
\begin{center}\begin{tikzpicture}
\matrix(m)[ampersand replacement=\&, matrix of math nodes, row sep=3em, column sep=2.5em, text height=1.5ex, text depth=0.25ex]
{P_1\&P_0\&X\&0\\}; 
\path[->]
(m-1-1) edge (m-1-2)
(m-1-2) edge (m-1-3)
(m-1-3) edge (m-1-4);
\end{tikzpicture}\end{center}   there is an exact sequence 
\begin{center}\begin{tikzpicture}
\matrix(m)[ampersand replacement=\&, matrix of math nodes, row sep=3em, column sep=2.5em, text height=1.5ex, text depth=0.25ex]
{0\&(X,\blank)\&(P_0,\blank)\&(P_1,\blank)\\}; 
\path[->]
(m-1-1) edge (m-1-2)
(m-1-2) edge (m-1-3)
(m-1-3) edge (m-1-4);
\end{tikzpicture}\end{center}
Since $P_i$ are finitely generated projective, it follows that $(P_i,\blank)\cong (P_i)^*\otimes\blank$.  Hence there is an exact sequence 
\begin{center}\begin{tikzpicture}
\matrix(m)[ampersand replacement=\&, matrix of math nodes, row sep=3em, column sep=2.5em, text height=1.5ex, text depth=0.25ex]
{0\&(X,\blank)\& P_0^*\otimes\blank\&P_1^*\otimes\blank\&\tr(X)\otimes\blank\&0\\}; 
\path[->]
(m-1-1) edge (m-1-2)
(m-1-2) edge (m-1-3)
(m-1-3) edge (m-1-4)
(m-1-4) edge (m-1-5)
(m-1-5) edge (m-1-6);
\end{tikzpicture}\end{center}
which is an injective resolution of $(X,\blank)$.  Since $D:=L^0\Y_R$ we have that $D(X,\blank)$ is determined by the top row in the following commutative diagram with exact rows

\dia{(P_1^*,\blank)\ar[r]\ar[d]^{\cong}&(P_0^*,\blank)\ar[r]\ar[d]^{\cong}&D(X,\blank)\ar[r]\ar[d]^\varphi&0\\
\blank\otimes P_1\ar[r]&\blank\otimes P_0\ar[r]&\blank\otimes X\ar[r]&0\\} From this, it follows that $\varphi$ is an isomorphism and hence$$D_R(X,\blank)\cong \blank\otimes X\qquad \qed$$ \end{pf}

\begin{lem}For any $X\in \mod(R^{op})$, $D_R(\blank\otimes X)\cong (X,\blank)$ where the finitely presented functor $D_R(\blank\otimes X)$ is being viewed as an object in $\tfp(\Mod(R),\ab)$\end{lem}

\begin{pf}Given $X\in \mod(R^{op})$, take presentation $P_1\to P_0\to X\to 0$.  There is a commutative diagram with exact rows \dia{\blank\otimes P_1\ar[r]\ar[d]^{\cong}& \blank\otimes P_0\ar[r]\ar[d]^{\cong} &\blank\otimes X\ar[r]\ar[d] ^{\cong}&0\\
(P_1^*,\blank)\ar[r]& (P_0^*,\blank)\ar[r]& \blank\otimes X\ar[r]& 0} where left two vertical morphisms are isomorphisms because the $P_i$ are finitely generated projectives. Applying $D_R$ to the bottom row yields the following commutative diagram with exact rows, again using the fact that the $P_i$ are finitely generated projectives
\dia{0\ar[r]& D_R(\blank\otimes X)\ar[d]^{\cong}\ar[r]& D_R(P_0^*,\blank)\ar[d]^{\cong}\ar[r]& D_R(P_1^*,\blank)\ar[d]^{\cong}\\ 
0\ar[r]&D_R(\blank\otimes X)\ar[r]\ar[d]^\varphi& P_0^*\otimes\blank\ar[r]\ar[d]^{\cong}& P_1^*\otimes\blank\ar[d]^{\cong}\\
0\ar[r]& (X,\blank)\ar[r]& (P_0,\blank)\ar[r]& (P_1,\blank)}   Therefore $\varphi$ is an isomorphism and we obtain $$D_R(\blank\otimes X)\cong (X,\blank) \qquad  \qed$$\end{pf}

An immediate result of these two lemmas is one of two main reasons that we have shifted discussion to the category $\tfp(\Mod(R),\ab)$.  While $D_R$ is not a duality, its restriction to $\tfp(\Mod(R),\ab)$ is a duality.

\begin{thm}The exact functor $$D_R\:\fp(\Mod(R),\ab)\to \fp(\Mod(R^{op}),\ab)$$ restricts to an exact duality $$D_A\:\tfp(\Mod(R),\ab)\to \tfp(\Mod(R^{op}),\ab)$$



\end{thm}

\begin{pf}Since $D_R$ is contravariant and exact, $D_R^2$ is covariant and exact.  Take any $F\in \tfp(\Mod(R),\ab)$ and choose $X,Y\in \mod(R)$ and presentation $$(Y,\blank)\To (X,\blank)\To F\To 0$$  Applying $D^2_R$ and observing that on representable functors $D^2_R$ is isomorphic to the identity functor, we obtain a commutative diagram with exact rows
\dia{D^2_R(Y,\blank)\ar[d]_{\cong}\ar[r]& D^2_R(X,\blank)\ar[r]\ar[d]_{\cong}& D^2_R(F)\ar[r]\ar[d]^\varphi &0\\
(Y,\blank)\ar[r]& (X,\blank)\ar[r] &F\ar[r]& 0}  It follows that $\varphi$ is an isomorphism and hence $D_R^2F\cong F$.   Therefore $D_A$ which is the restriction of $D_R$ to $\tfp(\Mod(R),\ab)$ is a duality. $\qed$.  \end{pf}
We now have the following statement first discovered by Auslander:
\begin{lem}[Auslander, \cite{isosing}]Let $F\in \tfp(\Mod(R),\ab)$.  Then there exists $X,Y,Z\in \mod(R^{op})$ and injective resolution $$0\To F\To \blank\otimes X\To\blank \otimes Y\To\blank\otimes Z\To 0$$\end{lem}

\begin{pf}Since $F\in \tfp(\Mod(R),\ab)$, the dual functor $D_AF\in \tfp(\Mod(R^{op}),\ab)$.  Therefore, there exists $X,Y,Z\in \mod(R^{op})$ with presentation $$0\To (Z,\blank)\To (Y,\blank)\To (X,\blank)\To D_AF\To 0$$Applying the exact duality $D_A$ and using the isomorphisms from above we have the following commutative diagram with exact rows


\begin{center}\begin{tikzpicture}
\matrix(m)[ampersand replacement=\&, matrix of math nodes, row sep=3em, column sep=2.5em, text height=1.5ex, text depth=0.25ex]
{0\&D_A^2F\&D_A(X,\blank)\&D_A(Y,\blank)\&D_A(Z,\blank)\&0\\
0\&F\&\blank\otimes X\&\blank \otimes Y\&\blank\otimes Z\&0\\}; 
\path[->]
(m-1-1) edge node[auto]{}(m-1-2)
(m-1-2) edge node[auto]{}(m-1-3)
(m-1-3) edge node[auto]{}(m-1-4)
(m-1-4) edge node[auto]{}(m-1-5)
(m-1-2) edge node[left]{$\cong$}(m-2-2)
(m-2-1) edge node[auto]{}(m-2-2)
(m-2-2) edge node[auto]{}(m-2-3)
(m-2-3) edge node[below]{}(m-2-4)
(m-2-4) edge node[auto]{}(m-2-5)
(m-1-3) edge node[auto]{$\cong$}(m-2-3)
(m-1-4) edge node[auto]{$\cong$}(m-2-4)
(m-1-5) edge node[auto]{$\cong$}(m-2-5)
(m-2-5)edge node[auto]{}(m-2-6)
(m-1-5)edge node[auto]{$\cong$}(m-1-6);
\end{tikzpicture}\end{center}  This completes the proof.  $\qed$ \end{pf}

We have arrived at the second major reason for shifting our focus to $\tfp(\Mod(R),\ab)$.  The satellite endofunctors on the larger category $(\Mod(R),\ab)$ can be restricted to endofunctors on the smaller category.  The fact that we are focusing on functors from $\Mod(R)$ instead of $\mod(R)$ is crucial.  The satellites of any functor $F\:\Mod(R)\to \ab$ exist because $\Mod(R)$ has enough injectives and projectives.  This technical difficulty could be avoided by using the more general definition of satellies given by Fisher-Palmquist and Newell; however, that approach is slightly more involved.

\begin{thm}For any $F\in \tfp(\Mod(R),\ab)$ and for all $n\ge 0$: 
\begin{enumerate}
\item $S^nF\in \tfp(\Mod(R),\ab)$
\item $S_nF\in \tfp(\Mod(R),\ab)$
\end{enumerate}
As a result, using the results of Fisher-Palmquist and Newell, $(S^n,S_n)$ form an adjoint pair on the functor category $\tfp(\Mod(R),\ab)$.
\end{thm}

\begin{pf}The case $n=0$ is trivial and all other cases follow readily from the case $n=1$.  

Assume we have $X,Y\in \mod(R)$ and presentation  \dia{(Y,\blank)\ar[r]&(X,\blank)\ar[r]&F\ar[r]&0} Apply the right exact functor $S^1$ to get the following commutative diagram with exact rows 
\dia{S^1(Y,\blank)\ar[r]\ar[d]^{\cong}& S^1(X,\blank)\ar[r]\ar[d]^\cong& S^1F\ar[d]^{\cong}\ar[r]& 0\\
\ext^1(Y,\blank)\ar[r]&\ext^1(X,\blank)\ar[r]&S^1F\ar[r]&0}   But since $X,Y\in \mod(R)$, both $\ext^1(X,\blank)$ and $\ext^1(Y,\blank)$ are in $\tfp(\Mod(R),\ab)$ and since this category is abelian, $S^1F$ is totally finitely presented.  This shows $(1)$.

Suppose that $F\in\tfp(\Mod(R),\ab)$.  Since $R$ is coherent, there are finitely presented left modules $X,Y,Z$ and an exact sequence $$0\to F\to \blank\otimes X\to \blank\otimes Y\to \blank\otimes Z\to 0$$  It is easily seen that $S_1$ preserves left exact sequences of functors with respect to evaluation.  Therefore, applying $S_1$ yields the following commutative diagram with exact rows  \dia{0\ar[r]& S_1F\ar[r]\ar[d]^{\cong}& S_1(\blank\otimes X)\ar[r]\ar[d]^{\cong}& S_1(\blank \otimes Y)\ar[d]^{\cong}\\
0\ar[r]&S_1F\ar[r]&\tor_1(\blank, X)\ar[r]&\tor_1(\blank, Y)} Since both $\tor_1(\blank,X)$ and $\tor_1(\blank, Y)$ are totally finitely presented, it follows that $S_1F$ is also totally finitely presented.   This establishes (2) and completes the proof.  $\qed$
\end{pf}

 \begin{lem}For any $M\in \mod(R^{op})$, $$D_A(\ext^1(M,\blank))\cong \tor_1(\blank,M)$$\end{lem}
 
 \begin{pf} The syzygy sequence $$0\to \u M\to P\to M\to 0$$ yields exact sequence $$ (P,\blank)\to (\u M,\blank)\to \ext^1(M,\blank)\to 0.$$ Applying the exact functor $D_A$ yields commutative diagram with exact rows 
 \dia{0\ar[r]& D_A(\ext^1(M,\blank)\ar[r]\ar[d]^\varphi& D_A(\u M, \blank)\ar[r]\ar[d]^{\cong}& D_A(P,\blank)\ar[d]^{\cong}\\
 0\ar[r]& \tor_1(\blank,M)\ar[r]& \blank\otimes \u M\ar[r]& \blank\otimes P} it follows that $\varphi$ is an isomorphism and hence $D_A(\ext^1(M,\blank))\cong \tor_1(\blank,M)$.$\qed$\end{pf}

\begin{lem}For any $M\in \mod(R^{op})$, $$D_A(\tor_1(\blank,M))\cong \ext^1(M,\blank)$$\end{lem}

\begin{pf}
If $M\in \mod(R^{op})$, from the syzygy sequence $$0\to \u M\to P\to M\to 0$$ there is an exact sequence $$0\to \tor_1(\blank,M)\to \blank\otimes \u M\to \blank\otimes P$$  Applying the exact functor $D_A$ yields the following commutative diagram with exact rows \dia{D_A(\blank\otimes P)\ar[r]\ar[d]^{\cong}& D_A(\blank\otimes\u X)\ar[r]\ar[d]^{\cong}& D_A(\tor_1(\blank,M))\ar[r]\ar[d]& 0\\
(P,\blank)\ar[r]& (\u M,\blank)\ar[r]&\ext^1(M,\blank)\ar[r]& 0} it follows that $$D_A(\tor_1(\blank,M))\cong \ext^1(M,\blank)\qquad\qed$$\end{pf}

\section{The Auslander-Gruson-Jenson Duality}

The functor $D_A$ is actually a well known duality.  It was first discovered by Auslander and appears in \cite{isosing} and independently discovered by Gruson and Jensen and appears in \cite{gj}.  For this reason it is often referred to as the Auslander-Gruson-Jensen duality.  It turns out that this duality also makes an appearance in model theory.  This occurs after a sequence of individual contributions made by Prest, Herzog, and Burke.  We will not go into details here but the interested reader is referred to \cite{psl} for an overview or \cite{prestmt},  \cite{herzogduality}, and \cite{burke} for a full account of how the duality appears in model theory.

  \begin{prop}For all $n\ge 0$, 
  \begin{enumerate}
  \item $D_AS^n\cong S_nD_A$
  \item $D_AS_n\cong S^nD_A$
  \end{enumerate}
  \end{prop}
 \begin{pf}From the exact sequence $$(Y,\blank)\to (X,\blank)\to F\to 0$$ we apply the right exact functor $S^1$ to get the following commutative diagram with exact rows \dia{S^1(Y,\blank)\ar[r]\ar[d]^{\cong}& S^1(X,\blank)\ar[r]\ar[d]^{\cong}& S^1F\ar[r]\ar[d]^{\cong}& 0\\
  \ext^1(Y,\blank)\ar[r]& \ext^1(X,\blank)\ar[r]& S^1F\ar[r]& 0}  Applying $D_A$ yields the following commutative diagram with exact rows 
  \dia{0\ar[r]&D_AS^1F\ar[r]\ar[d]^{\cong}&D_A(\ext^1(X,\blank))\ar[r]\ar[d]^{\cong}&D_A(\ext^1(Y,\blank))\ar[d]^{\cong}\\
  0\ar[r]& D_AS^1F\ar[r]\ar[d]^{\cong} &\tor_1(\blank, X)\ar[r]\ar[d]^{\cong}& \tor_1(\blank,Y)\ar[d]^{\cong}\\
0\ar[r]& D_AS^1F\ar[r]& S_1(\blank\otimes X)\ar[r]& S_1(\blank\otimes Y)}  The bottom row of this diagram can also be obtained by applying $D_A$ to the exact sequence $$(Y,\blank)\to (X,\blank)\to F\to 0$$ which results in $$0\To D_AF\To \blank\otimes X\To \blank \otimes Y$$ and then applying the left exact functor $S_1$ $$0\to S_1D_AF\to S_1(\blank\otimes X)\to S_1(\blank\otimes Y)$$ Hence $$S_1D_AF\cong D_AS^1F$$ One may similarly show that $D_AS_1\cong S^1D_A$. The result will now follow by induction on $n$.  $\qed$\end{pf}

At this point we will now abandon using the notation $D_A$ and simply refer to this functor as $D$.  Notice that $S^n(M,\blank)\cong \ext^n(M,\blank)$.  Hence $$D\ext^n(M,\blank)\cong DS^n(M,\blank)\cong S_nD(M,\blank)\cong S_n(\blank\otimes M)\cong \tor_n(\blank,M)$$  Similarly, $D\tor_n(\blank,M)\cong \ext^n(M,\blank)$.  Many of the calculations involving $D$ from above can be found in both \cite{isosing} and \cite{cohhart}.  The main difference there is that both Hartshorne and Auslander use the category $\fp(\mod(R),\ab)$ instead of the category $\tfp(\Mod(R),\ab)$.   We reiterate that the reason for using the functors capable of dealing with big modules is the desire to understand the satellites of finitely presented functors in the most efficient way possible.  As a result, we are able to state and prove precisely the anticommutative relationship between $D$ and the satellites $S^n,S_n$.   These calculations cannot be found in the literature and are crucial to understanding how to extend linkage to the category of finitely presented functors.   We summarize all of the results in the following two theorems.

\begin{thm}For any coherent ring $R$, the functor $$D\:\tfp(\Mod(R^{op}),\ab))\to \tfp(\Mod(R),\ab)$$ is a duality satisfying the following properties for all $n\ge 0$:
\begin{enumerate}
\item $D(\ext^n(M,\blank))\cong \tor_n(\blank, M)$
\item $D(\tor_n(\blank,M))\cong \ext^n(M,\blank)$
\item $DS^n\cong S_nD$
\item $DS_n\cong S^nD$
\end{enumerate}\end{thm}

In order to prove the next theorem we need the following result.  
\begin{prop}[Auslander,\cite{coh}]If $H\:\Mod(R)\To \ab$ is any right exact functor, then the functor $\nat(\blank,H)\:\fp(\Mod(R),\ab)\To \ab$ is exact.\end{prop}

\begin{thm}The duality $D\:\tfp(\Mod(R),\ab)\to \tfp(\Mod(R^{op}),\ab)$ satisfies the following properties:
\begin{enumerate}
\item $DF(A)\cong \nat(F,\blank\otimes A)$
\item Given $F\in \tfp(\Mod(R),\ab)$, take presentation $(Y,\blank)\to (X,\blank)\to F\to 0$.  Then $DF$ is completely determined by the exact sequence $0\to DF\to X\otimes \blank\to Y\otimes \blank$.  
  \end{enumerate} As a result, 
\begin{enumerate}
\item $D$ is the duality first discovered by Auslander appearing in \cite{isosing} and independently discover by Gruson and Jensen appearing in \cite{gj}.
\item $D$ is the duality studied by Hartshorne in \cite{cohhart}.
\item $D$ is the duality appearing in model theory whose existence is due to the sequence of results in the field:
\begin{enumerate}
\item Prest defines $D$ of a pp-formula.  \cite{prestmt}
\item Herzog defines $D$ on the category of pp-pairs which he shows is an abelian category. \cite{herzogduality}
\item Burke shows that the the category of pp-pairs is equivalent to the category of finitely presented functors and that $D$ is the Auslander-Gruson-Jensen duality.  \cite{burke}
\end{enumerate}
\end{enumerate}\end{thm}

\begin{pf}For any $F\in \tfp(\Mod(R),\ab)$ take presentation $$(Y,\blank)\to (X,\blank)\to F\to 0$$ and apply the exact functor $\nat(\blank, \blank\otimes A)$ to get exact sequence $$0\to \nat(F,\blank\otimes A)\to \nat\big((X,\blank),\blank\otimes A\big)\to \nat\big((Y,\blank),\blank\otimes A\big)$$ which by the Yoneda lemma is equivalent to  $$0\to (F,\blank\otimes A)\to X\otimes A\to Y\otimes A$$  Finally, applying the exact functor $D$ to the same presentation of $F$ yields exact sequence $$0\to DF\to D(X,\blank)\to D(Y,\blank)$$ or equivalently $$0\to DF\to X\otimes \blank\to Y\otimes \blank.$$  Evaluating at $A$ yields exact sequence $$0\to DF(A)\to X\otimes A\to Y\otimes A.$$  Hence $DF(A)\cong \nat(F,\blank\otimes A)$ as claimed.  $\qed$  \end{pf}

It should be pointed out that although $$D\:\tfp(\Mod(R),\ab)\to \tfp(\Mod(R^{op}),\ab)$$ is equivalent to the Auslander-Gruson-Jensen duality, the original functor $D_R\:\fp(\Mod(R),\ab)\to \fp(\Mod(R^{op}),\ab)$ does not seem to obey the formula $D_RF(A)=\nat(F,A\otimes \blank)$.  Moreover, we recovered the Auslander-Gruson-Jensen duality by exploiting the fact that $\fp(\Mod(R),\ab)$ has enough injectives which avoids any mention of the tensor product.   Even though these two larger functors do not appear to be in general isomorphic, we do have the following relationship between the two:

\begin{thm}Let $D_{aus}$ denote the functor defined by Auslander by the formula $DF(A)=\nat(F,\blank\otimes A)$.  There is a natural transformation $$D_{aus}\To D_R$$\end{thm}

\begin{pf}Take any injective resolution any $F\in \fp(\Mod(R),\ab)$ $$0\To F\To I^0\To I^1$$  For any $A\in\Mod(R^{op})$, this gives rise to the diagram of abelian groups

\begin{center}\begin{tikzpicture}
\matrix(m)[ampersand replacement=\&, matrix of math nodes, row sep=3em, column sep=2.5em, text height=1.5ex, text depth=0.25ex]
{\nat(I^1,\blank\otimes A)\&\nat(I^0,\blank\otimes A)\&D_{aus}F(A)\&0\\
(I^1(R),A)\&(I^0(R),A)\&D_RF(A)\&0\\}; 
\path[->]
(m-1-1)edge node[left]{$\ev_R$}(m-2-1)
(m-1-1) edge node[auto]{}(m-1-2)
(m-1-2) edge node[auto]{}(m-1-3)
(m-1-3) edge node[auto]{}(m-1-4)
(m-1-2) edge node[left]{$\ev_R$}(m-2-2)
(m-2-1) edge node[auto]{}(m-2-2)
(m-2-2) edge node[auto]{}(m-2-3)
(m-2-3) edge node[below]{}(m-2-4)
(m-1-3) edge node[right]{$\exists !\varphi$}(m-2-3);
\end{tikzpicture}\end{center}

which is easily seen to be natural in $A$.  $\qed$

\end{pf}

%

We end this section by showing how $D$ may be used to calculate the defect $w$ of a totally finitely presented functor.   The defect is the contravariant exact functor $w\:\fp(\Mod(R),\ab)\To \Mod(R)$ satisfying $w(X,\blank)=X$.  This determines $w$ completely.  Again, recall that we are assuming that $R$ is coherent.  We begin by calculating $w(\blank\otimes X)$.  First, take presentation $$P_1\to P_0\to X\to 0.$$  This yields exact sequence $$0\to X^*\to P_0^*\to P_1^*\to \tr(X)\to 0$$  Embedding this into $\fp(\mod(R),\ab)$ yields exact sequence $$0\to (\tr(X),\blank)\to (P_1^*,\blank)\to (P_0^*,\blank)\to F\to 0.$$  Since the $P_i$ are finitely generated projective, this is equivalent to the following exact sequence $$0\to (\tr(X),\blank)\to \blank\otimes P_1\to \blank\otimes P_0\to F\to 0.$$  Since we have exact sequence $$\blank\otimes P_1\to \blank\otimes P_0\to \blank\otimes X\to 0$$ it follows that $F\cong\blank\otimes X$ and we have exact sequence $$0\to (\tr(X),\blank)\to (P_1^*,\blank)\to (P_0^*,\blank)\to \blank\otimes X\to 0.$$  By applying $w$ we get exact sequence $$0\to w(\blank\otimes X)\to P_0^*\to P_1^*\to \tr(X)\to 0$$ from which it follows that $w(\blank\otimes X)=X^*$.  

With this information, for any finitely presented functor $F$ take presentation $$(Y,\blank)\to (X,\blank)\to F\to 0$$ and apply $D$ yielding exact sequence $$0\to D(F)\to \blank\otimes X\to \blank\otimes Y.$$  Applying $w$ commutative diagram with exact rows:\begin{center}\begin{tikzpicture}
\matrix(m)[ampersand replacement=\&, matrix of math nodes, row sep=3em, column sep=2.5em, text height=1.5ex, text depth=0.25ex]
{Y^*\&X^*\&w\big[D(F)\big]\&0\\
(Y,R)\&(X,R)\&F(R)\&0\\}; 
\path[->]
(m-1-1)edge node[left]{$1$}(m-2-1)
(m-1-1) edge node[auto]{}(m-1-2)
(m-1-2) edge node[auto]{}(m-1-3)
(m-1-3) edge node[auto]{}(m-1-4)
(m-1-2) edge node[left]{$1$}(m-2-2)
(m-2-1) edge node[auto]{}(m-2-2)
(m-2-2) edge node[auto]{}(m-2-3)
(m-2-3) edge node[below]{}(m-2-4)
(m-1-3) edge node[auto]{$\cong$}(m-2-3);
\end{tikzpicture}\end{center} Hence $w\big[D(F)\big]\cong F(R)$.  Now using the fact that $D^2F\cong F$, one has
$$w(F)\cong w(D^2(F))
\cong w\big[D(DF)\big]
\cong DF(R)
$$ Hence the defect of a totally finitely presented functor is completely determined by the dual of the functor evaluated at the ring. 
 
\begin{thm}For any $F\in \tfp(\Mod(R),\ab)$ 
\begin{enumerate}
\item $w\big[D(F)\big]\cong F(R)$
\item $w(F)\cong DF(R)$
\end{enumerate}
\end{thm}

\begin{cor}For any $F\in \tfp(\Mod(R),\ab)$ and for all $n\ge 0$,  $$w(S_nF)\cong S^n(DF)(R)$$\end{cor}

\begin{pf}This follows from the fact that $S^nD\cong D S_n$ and the preceding theorem.  As a result $$w(S_n F)=D(S_nF)(R)
=S^n(DF)(R)$$ as claimed.  $\qed$\end{pf}

\begin{cor}For all $n\ge 0$ and for all $M\in \mod(R^{op})$, \begin{eqnarray*}w\big(\tor_n(\blank,M)\big)&\cong&\ext^n(M,R)\end{eqnarray*}\end{cor}

\section{Linkage of Finitely Presented Functors}

The functor $D$ together with the satellites will allow us to extend the definition of linkage to the category of totally finitely presented functors in a suitable way.  Many of the techniques of this section are established in \cite{stab}.

Recall  that for any module $M$, $S^1\ext^1(M,\blank)\cong \ext^1(\u M,\blank)$.  This is why we suspect that $S^1$ should play the role of $\u$ in the functor category.  Given $M\in \mod(R)$, we have $D\ext^1(M,\blank)=\tor_1(M,\blank)$.  This means that $D$ does not play the role of the transpose; however, we now recall
\begin{prop}[Auslander-Bridger, \cite{stab}]Let $M\in \mod(R)$, then $$S^1\tor_1(M,\blank)\cong \ext^1(\tr M,\blank)$$ \end{prop}

\begin{prop}The role of the transpose at the level of totally finitely presented functors is played by the functor $$S^1D\:\tfp(\Mod(R),\ab)\To \tfp(\Mod(R^{op}),\ab)$$  In particular $$S^1D\big(\ext^1(M,\blank)\big)\cong \ext^1(\tr M,\blank)$$
\end{prop}

 \begin{thm}The diagram of functors \begin{center}\begin{tikzpicture}[description/.style={fill=white,inner sep=2pt}]  
\matrix (m) [ampersand replacement= \&,matrix of math nodes, row sep=3em, 
column sep=6em,text height=1.5ex,text depth=0.25ex] 
{\underline{\mod}(R)\&\underline{\mod}(R^{op})\\
\&\\
\tfp(\Mod(R),\ab)\&\tfp(\Mod(R^{op}),\ab)\\}; 
\path[->,thick,blue, font=\scriptsize]
(m-1-1) edge[bend left] node[auto]{$\tr$}(m-1-2)
(m-1-2) edge[bend left] node[auto]{$\tr$}(m-1-1);
\path[->,thick, red, font=\scriptsize,  loop/.style={min distance=2cm,in=45, out=135, looseness=13}]
(m-1-2) edge[loop] node[description]{$\Omega$}(m-1-2);
\path[->,thick, red, font=\scriptsize,  loop/.style={min distance=2cm,in=135, out=45, looseness=13}]
(m-1-1) edge[loop] node[description]{$\Omega$}(m-1-1);
\path[->, dashed, thick,blue, font=\scriptsize]
(m-3-1) edge[bend left] node[auto]{$S^1D$}(m-3-2)
(m-3-2) edge[bend left] node[auto]{$S^1D$}(m-3-1);
\path[->, dashed, thick, red, font=\scriptsize,  loop/.style={min distance=2cm,in=-45, out=-135, looseness=13}]
(m-3-2) edge[loop] node[description]{$S^1$}(m-3-2);
\path[->, dashed,thick, red, font=\scriptsize,  loop/.style={min distance=2cm,in=-135, out=-45, looseness=13}]
(m-3-1) edge[loop] node[description]{$S^1$}(m-3-1);
\path[->,dashed, thick, font=\scriptsize]
(m-1-1) edge node[left]{$\hr$} (m-3-1);
\path[->,dashed, thick, font=\scriptsize]
(m-1-2) edge node[auto]{$\hr$} (m-3-2);
\end{tikzpicture}\end{center} satisfies the following properties:

\begin{enumerate}
\item The vertical arrows are contravariant embeddings.
\item The following commutativity relations are satisfied:
\begin{enumerate}
\item $S^1D\hr=\hr\tr$
\item $S^kS^1D\tr=\hr \u^k\tr$
\end{enumerate}
\item There is a counit of adjunction $S^2S_2\to 1$.
\end{enumerate} \end{thm}

Note $S^1D$ is the analog of $\tr$ and $S^1$ is the analog of $\u$.  As such, one may define a functor $F\in \tfp(\Mod(R),\ab)$ to be horizontally linked if $F\cong S^1 S^1D S^1S^1D(F)$.  This is equivalent to saying $F$ is linked if $F\cong S^2DS^2D(F)$.  Since $S^nD\cong DS_n$,  $F$ is linked if $F\cong S^2S_2D^2(F)\cong S^2S_2(F)$.  Moreover, this is equivalent to requiring the counit of adjunction of the adjoint pair $(S^2,S_2)$ is an isomorphism at $F$.  In other words, one may give the following:

\begin{defn}A functor $F\in \tfp(\Mod(R),\ab)$ is \dbf{horizontally linked} if the counit of adjunction \dia{S^2S_2F\ar[r]^{\quad c_F}&F}  evaluated at $F$ is an isomorphism. \end{defn}

As an immediate result, we see that the only functors that have a chance of being linked are the injectively stable functors.  One of the first requirements we imposed was that the new definition of linkage be consistent with that given for the stable category $\underline{\mod}(R)$.  We in fact have:

\begin{thm}A module $M\in \underline{\mod}(R)$ is linked if and only if $\ext^1(M,\blank)$ is linked.\end{thm}

\begin{pf}Both occur if and only if $\u\tr\u\tr M\cong M$.  $\qed$\end{pf}

In \cite{linkage}, Martsinkovsky and Strooker establish that all modules of $G$-dimension zero are horizontally linked.  The notion of $G$-dimension was introduced by Auslander and Bridger in \cite{stab} and the definition is given for half exact functors.  Suppose that $\A$ is an abelian category with enough projectives and injectives.  A half exact functor $F\:\A\to \B$ is said to have $G$-dimension zero if all of its satellites are both projectively and injectively stable.  A module is said to have $G$-dimension zero if the functor $\underline{\hom}(M,\blank)$ has $G$-dimension zero as a functor.  That modules of $G$-dimension zero are linked is an immediate result of the more general statement that we can now make:

\begin{thm}All half exact finitely presented functors of $G$-dimension zero are horizontally linked.\end{thm}

\begin{pf}If $F$ has $G$-dimension zero, then both $S_1F$ and $F$ are injectively stable.  Therefore \begin{eqnarray*}S^2S_2F&\cong&S^1 S^1S_1 (S_1F)\\
&\cong&S^1\overline{S_1 F}\\
&\cong& S^1(S_1 F)\\
&\cong& \overline{F}\\
&\cong& F\end{eqnarray*}\end{pf}

We end this section with a classification of horizontally linked half exact functors.  It turns out that these all arise as extension functors.  This means that if there are horizontally linked functors that are not isomorphic to extension functors, then these functors must not be half exact; however, it is not clear that even the extension functors that are horizontally linked arise from finitely presented modules and hence there are open questions here.   

In \cite{freyd}, Freyd shows that if $\A$ is an abelian category closed under denumerable sums, then any direct summand of $\ext^1(M,\blank)$ is itself an extension functor.  That is if there exists a section $F\To \ext^1(M,\blank)$, then $F\cong \ext^1(B,\blank)$ for some $Y\in \A$.   In \cite{coh}, Auslander shows that any half exact functor $F\in \fp(\A,\ab)$ for which $w(F)=0$ is a direct summand of an extension functor $\ext^1(M,\blank)$.  As a result Auslander establishes the following:

\begin{prop}[Auslander, \cite{coh}]Suppose that $\A$ is an abelian category closed under denumerable sums. Any half exact finitely presented functor $F\in \fp(\A,\ab)$ for which $w(F)=0$ is an extension functor.\end{prop}

\begin{lem}Suppose that $G$ is a half exact functor in $\fp(\Mod(R),\ab)$.  Then $S^1G$ is an extension functor.\end{lem}

\begin{pf}Since $G$ is half exact, $S^1G$ is half exact.  Since $S^1G$ is injectively stable, $w(S^1G)=0$. Therefore $S^1G$ is an extension functor. $\qed$\end{pf}

\begin{thm}If $F\in\tfp(\Mod(R),\ab)$ is half exact and horizontally linked, then $F$ is an extension functor.\end{thm}

\begin{pf}Since $F$ is half exact, so is $S^1S_2F$.  Since  $S^2S_2F=S^1(S^1S_2F)$, it follows from the preceding lemma that $S^2S_2F\cong \ext^1(X,\blank)$.  As a result $F\cong \ext^1(X,\blank)$.  $\qed$\end{pf}

This last result raises an interesting question.  Are there large modules $M\in \Mod(R)$ for which $\ext^1(M,\blank)$ is a horizontally linked totally finitely presented functor?  It is easy to see that if $M$ is finitely presented, then $\ext^1(M,\blank)$ is totally finitely presented and hence $\ext^1(M,\blank)$ will be horizontally linked if and only if $M$ is horizontally linked in $\underline{\mod}(R)$; however, this does not exclude the possibility that $\ext^1(M,\blank)$ may be totally finitely presented while $M$ is a large module.  In fact, we already have mentioned that if $\ext^1(M\oplus R^\alpha,\blank)\cong \ext^1(M,\blank)$ so the real question is whether $\ext^1(M,\blank)$ being finitely presented implies that $M\cong N\oplus P$ where $N$ is finitely presented and $P$ is projective.  We end this paper with a classification of half exact finitely presented $G$-dimension zero functors.

\begin{thm}A half exact finitely presented functor $F\in\tfp(\Mod(R),\ab)$ has $G$-dimension zero if and only if $F$ is an extension functor $\ext^1(M,\blank)$ for some module $M$ which has $G$-dimension zero.\end{thm}

\bibliographystyle{amsplain}
\bibliography{slayer}

\end{document}